\documentclass{gtart}

\def\ifplaintex{\expandafter\ifx\csname documentclass\endcsname\relax}

\def\gtp{{\mathsurround=0pt\it $\cal G\mskip-2mu$eometry \&\ 
$\cal T\!\!$opology $\cal P\!$ublications}}  % GT publications

\def\recd{{\small Received:\qua\receiveddate\ifx\reviseddate\relax
\else\qquad Revised:\qua\reviseddate\fi\par}} 

\def\lognumber#1{\def\thelognumber{#1}}
\def\volumenumber#1{\def\thevolumenumber{#1}}
\def\volumeyear#1{\def\thevolumeyear{#1}}
\def\papernumber#1{\def\thepapernumber{#1}}
\def\pagenumbers#1#2{\def\startpage{#1}\def\finishpage{#2}}
\def\published#1{\def\publishdate{#1}}

\def\received#1{\def\receiveddate{#1}}
\def\revised#1{\def\reviseddate{#1}}
\def\accepted#1{\def\accepteddate{#1}}
\def\asciititle#1{\def\theasciititle{#1}}
\def\covertitle#1{\def\thecovertitle{#1}}

\long\def\asciiabstract#1{\long\def\theasciiabstract{#1}}

\let\\\par\let\thelognumber\relax\let\thevolumenumber\relax
\let\thepapernumber\relax\let\thevolumeyear\relax\let\startpage\relax
\let\finishpage\relax\let\publishdate\relax\let\receiveddate\relax
\let\reviseddate\relax\let\accepteddate\relax\let\theasciititle\relax
\let\thecovertitle\relax\let\theasciiauthors\relax
\let\theasciiabstract\relax

\let\theasciiemail\relax

\ifplaintex
\font\logobig=cmssbx10 scaled 3836
\font\logomed=cmssbx10 scaled 2557
\else
\font\logobig=cmssbx10 scaled 4200
\font\logomed=cmssbx10 scaled 2800
\fi

\long\def\makeagttitle{   %%% start of definition of \makeagttitle
\count0=\startpage
\agt\hfill      %   Journal title (top left) 
\hbox to 45truept{\vbox to 0pt{\vglue -13truept{\logomed A\kern -.37em{\logobig 
T}\kern -.38em G}\vss}\hss}
\break
{\small Volume \thevolumenumber\ (\thevolumeyear)
\startpage--\finishpage\nl
Published: \publishdate}

\vglue .25truein

{\parskip=0pt\leftskip 0pt plus
1fil\def\\{\par\smallskip}{\Large\bf\thetitle}\par\medskip} \vglue
0.05truein

{\parskip=0pt\leftskip 0pt plus 1fil\def\\{\par}{\sc\theauthors}
\par\medskip}%
 
\vglue 0.03truein 

{\small\leftskip 25truept\rightskip 25truept{\bf Abstract}\stdspace\theabstract

{\bf AMS Classification}\stdspace\theprimaryclass
\ifx\thesecondaryclass\relax\else; \thesecondaryclass\fi\par
{\bf Keywords}\stdspace \thekeywords\par}\vglue 7truept

}   %%%% end of definition of \makeagttitle

\ifplaintex
\font\phead=cmsl9 scaled 950
\font\pnum=cmbx10 scaled 913
\font\pfoot=cmsl9 scaled 950
\headline{\vbox to 0pt{\vskip -4.5mm\line{\small\phead\ifnum
\count0=\startpage ISSN 1472-2739 (on-line) 1472-2747 (printed)
\hfill {\pnum\folio}\else\ifodd\count0\def\\{ }% 
\ifx\theshorttitle\relax\thetitle\else\theshorttitle\fi\hfill{\pnum\folio}
\else\def\\{ and }{\pnum\folio}\hfill\ifx\theshortauthors\relax\theauthors
\else\theshortauthors\fi\fi\fi}\vss}}
\footline{\vbox to 0pt{\vglue 0mm\line{\small\pfoot\ifnum\count0=\startpage
\copyright\ \gtp\hfill\else
\agt, Volume \thevolumenumber\ (\thevolumeyear)\hfill\fi}\vss}}
\else
\font=cmsl9 scaled 1050
\font\lnum=cmbx10 
\font=cmsl9 scaled 1050
\makeatletter
\def\@oddhead{{\small\lhead\ifnum\count0=\startpage ISSN 1472-2739 
(on-line) 1472-2747 (printed)\hfill {\lnum\number\count0}\else\ifodd\count0
\def\\{ }\ifx\theshorttitle\relax \thetitle \else\theshorttitle\fi\hfill
{\lnum\number\count0}\else\def\\{ and }{\lnum\number\count0}
\hfill\ifx\theshortauthors\relax 
\theauthors\else\theshortauthors\fi\fi\fi}}\def\@evenhead{\@oddhead}
\def\@oddfoot{\small\lfoot\ifnum\count0=\startpage\copyright\ \gtp\hfill\else
\agt, Volume \thevolumenumber\ (\thevolumeyear)\hfill\fi}
\def\@evenfoot{\@oddfoot}
\makeatother
\fi
\let\maketitlepage\makeagttitle
\let\makeshorttitle\maketitlepage
\let\maketitle\maketitlepage

\newwrite\gtoutfile
\long\gdef\makeheadfile{  %%% start of definition of \makeheadfile
{\def\\{, }\def\s{ }
\immediate\openout\gtoutfile head.xxx
\immediate\write\gtoutfile{To: math@arxiv.org}
\immediate\write\gtoutfile{Subject: put OR rep NNNNN:ppppp}
\immediate\write\gtoutfile{--text follows this line--}
\immediate\write\gtoutfile{Proxy-for: \ifx\theasciiauthors\relax
\theauthors\else\theasciiauthors\fi\s<\ifx\theasciiemail\relax\theemail\else\theasciiemail\fi>}
\immediate\write\gtoutfile{\noexpand\\}
\immediate\write\gtoutfile{Authors: \ifx\theasciiauthors\relax
\theauthors\else\theasciiauthors\fi}
{\def\\{ }\immediate\write\gtoutfile{Title: \ifx\theasciititle\relax
\thetitle\else\theasciititle\fi}}
\immediate\write\gtoutfile{Subj-class: GT or SG, GR etc}
\immediate\write\gtoutfile{MSC-class: \theprimaryclass\ifx\thesecondaryclass\relax\else, \thesecondaryclass\fi}
\immediate\write\gtoutfile{Journal-ref: Algebr. Geom. Topol. \thevolumenumber\s
(\thevolumeyear) \startpage-\finishpage}
\immediate\write\gtoutfile{Comments: Published by Algebraic and
Geometric Topology at}
\immediate\write\gtoutfile{\s\s\s  http://www.maths.warwick.ac.uk/agt/AGTVol\thevolumenumber/agt-\thevolumenumber-\thepapernumber.abs.html}
\immediate\write\gtoutfile{\noexpand\\}
\immediate\write\gtoutfile{}
\ifx\theasciiabstract\relax
\immediate\write\gtoutfile{\theabstract}\else
\immediate\write\gtoutfile{\theasciiabstract}\fi
\immediate\write\gtoutfile{}
\immediate\write\gtoutfile{\noexpand\\}
\immediate\write\gtoutfile{}
\immediate\closeout\gtoutfile}}  %%% end of definition of \makeheadfile

\def\maketitlepage{\makeagttitle\makeheadfile}
\let\makeshorttitle\maketitlepage
\let\maketitle\maketitlepage

\lognumber{21}
\volumenumber{3}
\volumeyear{2003}
\papernumber{21}
\published{22 June 2003}
\pagenumbers{593}{622}
\received{23 January 2003}
\revised{27 February 2003}
\accepted{23 April 2003}

\usepackage{amssymb, epsf}

\newtheorem{thm}{Theorem}[section] %the resolution could also be [subsection] 

\newtheorem{lemma}[thm]{Lemma} 
\newtheorem{claim}[thm]{Claim} 
\newtheorem{prop}[thm]{Proposition} 
\theoremstyle{remark}
 
\newtheorem{defn}[thm]{Definition} 
\newtheorem{rem}[thm]{Remark} 
\newcommand{\<}{\langle} 
\renewcommand{\>}{\rangle} 
 
\def\si{\sigma} 
\def\de{\delta} 
\def\De{\Delta} 

\def\C{{\mathbb C}} 
 
\def\cL{{\mathcal L}} 
\def\F{{\mathbb F}}

\def\N{{\mathbb N}} 
\def\P{{\mathbb P}}

\def\Z{{\mathbb Z}}

\def\({\left(} 
\def\){\right)} 
\begin{document} 
 
\title{Plane curves and their fundamental groups:\\Generalizations of Uluda\u g's construction} 
\asciititle{Plane curves and their fundamental groups:\\Generalizations of Uludag's construction} 
\covertitle{Plane curves and their fundamental groups:\\Generalizations 
of Uluda\noexpand\u g's construction} 
\shorttitle{Plane curves and their fundamental groups}
 
\author{David Garber}

\address{Institut Fourier, BP 74, 38402 Saint-Martin D'Heres CEDEX, FRANCE} 
\email{garber@mozart.ujf-grenoble.fr} 
  
\begin{abstract} 
In this paper we investigate Uluda\u g's method for constructing  
new curves whose fundamental groups are central extensions of 
the fundamental group of the original curve by finite cyclic groups. 
 
In the first part, we give some generalizations to his method in order 
to get new families of curves with controlled fundamental groups. 
In the second part, we discuss some properties of groups 
which are preserved by these methods. Afterwards, we describe 
precisely the families of curves which can be obtained by applying the  
generalized methods to several types of plane curves.    
We also give an application of the general methods for constructing
new Zariski pairs.
\end{abstract} 
\asciiabstract{ 
In this paper we investigate Uludag's method for constructing  
new curves whose fundamental groups are central extensions of 
the fundamental group of the original curve by finite cyclic groups. 
 
In the first part, we give some generalizations to his method in order 
to get new families of curves with controlled fundamental groups. 
In the second part, we discuss some properties of groups 
which are preserved by these methods. Afterwards, we describe 
precisely the families of curves which can be obtained by applying the  
generalized methods to several types of plane curves.    
We also give an application of the general methods for constructing
new Zariski pairs.}

 
\maketitle 
 
\section{Introduction} 

The fundamental group of complements of plane curves is a very
important topological invariant with many different applications.
This invariant was used by Chisini \cite{chisini}, Kulikov \cite{Kul}
and Kulikov-Teicher \cite{KuTe} in order
to distinguish between connected components of the moduli space of
surfaces of general type. Moreover, the Zariski-Lefschetz hyperplane
section theorem (see \cite{milnor}) showed that 
$$\pi_1 (\P^N - S) \cong \pi_1 (H - H \cap S)$$
where $S$ is an hypersurface and $H$ is a generic 2-plane. Since 
$H \cap S$ is a plane curve, this invariant can be used also for computing the
fundamental group of complements of hypersurfaces in $\P^N$. 

A different direction for the need of fundamental groups' computations is for
getting more examples of Zariski pairs \cite{Z1,Z2}. A pair of plane
curves is called {\it a Zariski pair} if they have the 
same combinatorics (i.e. the same singular points and the same arrangement of
irreducible components), but their complements are not homeomorphic.
Several families of Zariski pairs were presented by Artal-Bartolo
\cite{AB,AB-C}, Degtyarev \cite{Deg}, Oka \cite{oka} and
Shimada \cite{shimada1,shimada2,shimada3}. 
Tokunaga and his coauthors thoroughly investigated
Zariski pairs of curves of degree $6$ (see \cite{ACCT,A-T,Tok1,Tok2} and 
\cite{Tok3}). 
Some candidates for weak Zariski pairs (where we take into account only the
types of the singular points) can also be found in \cite{GTV}, 
where any pair of arrangements with the same signature but with
different lattices can serve as a candidate for a Zariski pair (it is only
needed to be checked that the pair of arrangements have non-isomorphic
fundamental groups).

It is also interesting to explore new finite non-abelian groups which
are serving as fundamental groups of complements of plane curves.

Uluda\u g \cite{uludag1,uludag2} presents a way to obtain 
new curves whose fundamental groups are central  extensions 
of the fundamental group of a given curve. Using his method, one
can produce a family of examples of Zariski pairs from a given
Zariski pair (see also Section \ref{zariski_pairs} here). 
His main result is: 
\begin{thm}[Uluda\u g] 
Let $C$ be a plane projective curve and $G=\pi _1 (\P ^2 -C)$. 
Then for any $n \in \N$, there is a curve $\tilde C \subset \P ^2$ 
birational to $C$ such that  $\tilde G=\pi _1 (\P ^2 -\tilde C)$ 
is a central extension of $G$ by $\Z/(n+1)\Z$: 
$$1 \rightarrow \Z/(n+1)\Z \rightarrow \tilde G \rightarrow G \rightarrow 1$$  
In particular, if $C$ is irreducible so is $\tilde C$ as well. 
\end{thm} 
   
A natural question is which curves and fundamental
groups can be obtained by this method. Also, one might ask if this
method can be generalized, and what will be the effect of the
general method on the relation between the fundamental groups
of the original curve and the resulting curve.

In this paper we first generalize Uluda\u g's method to 
get new families of curves whose fundamental groups are
controlled by the fundamental group of the original curve in the same manner. 
Precisely, instead of using only 
two fibers for performing the elementary transformations between Hirzebruch 
surfaces $F_n$, we allow any finite number of different fibers. 
Afterwards, we list properties of
groups which are preserved by the methods. Also, we describe the curves
obtained by the application of these methods to several types of plane curves.  
Then we present some infinite families of new Zariski pairs which
can be obtained by the application of these methods.

Among the interesting results in this paper is the exploring of 
families of curves with deep singularities which yet have 
cyclic groups as the fundamental group of the complement (see the
beginning of Section \ref{properties}). Also,  
using the general construction, one can obtain  
more curves with finite and non-abelian fundamental groups 
and more Zariski pairs (see Section \ref{zariski_pairs}).

The paper is organized as follows.  
In Section \ref{uludag} we present Uluda\u g's original construction.  
In Section \ref{main}, we present  
some generalizations of Uluda\u g's construction, and we prove that also 
in the general constructions, the obtained curve has a fundamental 
group which is a central extension of the original curve's fundamental 
group by a finite cyclic group.   
Section \ref{properties} deals with properties of groups 
which are preserved while applying the constructions.  
In Section \ref{extension_cyclic}, we describe precisely the families of  
plane curves which can be obtained by the general constructions 
and then we calculate the degrees of the new curves.
At the end of this section, we describe some specific families of
curves obtained by applying the constructions to several types of
plane curves. In Section \ref{zariski_pairs} we present some  
new examples of Zariski pairs obtained by these constructions.    
 
\section{Uluda\u g's method}\label{uludag} 
The idea of the method is the following (it was partially  
 introduced by Artal-Bartolo \cite{AB1} and Degtyarev \cite{Deg}, 
and the  sequel was developed by Uluda\u g \cite{uludag1,uludag2}). 
If a curve $C_2$ is obtained from a curve $C_1$ by means of a Cremona  
transformation $\psi: \P^2 \to \P^2$, then $\psi$ induces an 
isomorphism 
$$\P^2 - (C_2 \cup A) \cong \P^2 - (C_1 \cup B)$$ 
where $A$ and $B$ are certain curve arrangements. Hence there is an induced 
isomorphism between their fundamental groups. The fundamental 
groups of the curves themselves are easy to compute by adding the 
relations which correspond to the arrangements. 
 
Now, if we start with a curve  $C_1$ whose fundamental group is known, 
one can find a curve $C_2$ whose fundamental group has not yet been known as 
being a fundamental group of a plane curve.

For these Cremona transformations one can use Hirzebruch surfaces 
$F_n$, $n \in \N \cup \{0\}$. 
In principle, the Hirzebruch surfaces are $\P^1$-bundles over 
$\P^1$. It is known that two Hirzebruch surfaces can be distinguished 
by the self-intersection of their exceptional section (for $F_n$ the 
self-intersection of its exceptional section is $-n$). 
 
\begin{figure}[ht!]
\epsfysize=8cm 
\centerline{\epsfbox{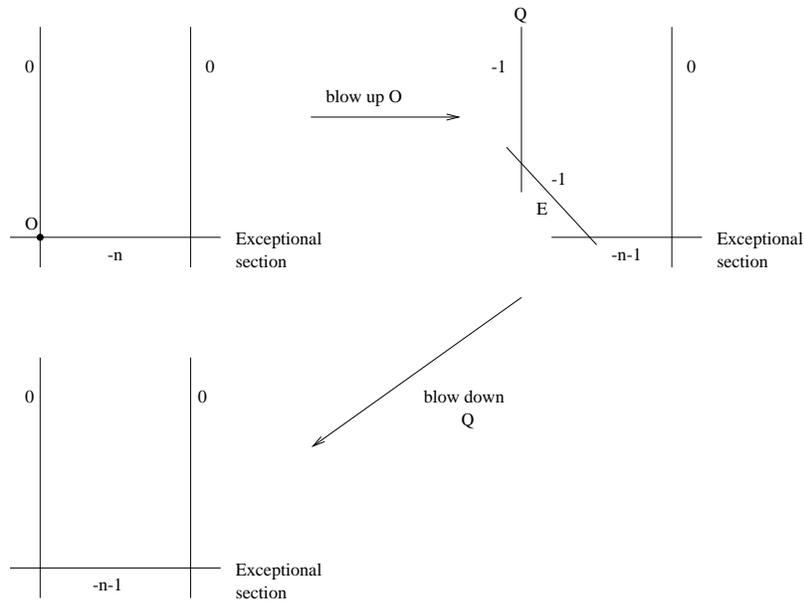}} 
\caption{Elementary transformation from $F_n$ to $F_{n+1}$}\label{fig1} 
\end{figure}

\begin{figure}[ht!]
\epsfysize=8cm 
\centerline{\epsfbox{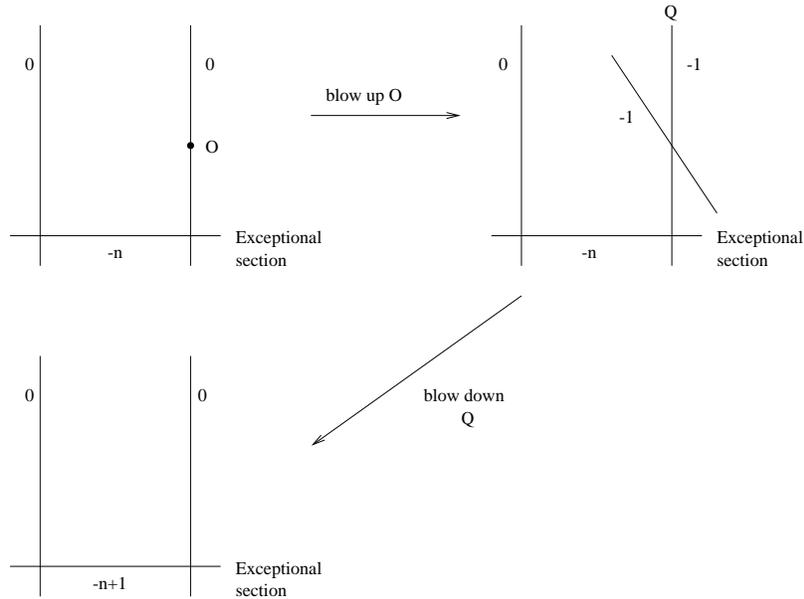}} 
\caption{Elementary transformation from $F_n$ to $F_{n-1}$}\label{fig2} 
\end{figure}

There are two types of elementary transformations, one transforms $F_n$ 
to $F_{n+1}$ for all positive $n$, and the other transforms $F_n$ 
to $F_{n-1}$ for all positive $n$ (these transformations are also called 
{\it Nagata elementary transformations}, see \cite{Nag}). 
The first transformation blows up a point 
$O$ on the exceptional section, and then blows down the proper 
transform of the fiber passed through $O$.   
The second transformation blows up a 
point $O$ on one of the fibers $Q$, outside the exceptional section,   
and then blows down the proper transform of the fiber $Q$. 
The two transformations are schematically presented in 
Figures \ref{fig1} and \ref{fig2} respectively. 
 
Uluda\u g used a special type of Cremona transformations which can be 
described as follows. We start with a 
curve $C$ in $\P^2$, and an additional line $Q$ which intersects the curve  
transversally. We choose another line $P$ which also intersects the curve  
transversally, and meets $Q$ outside  the curve $C$.  
Then we blow up the intersection point of the two lines. 
This yields the Hirzebruch surface $F_1$. Then we apply $n$ 
elementary transformations of the first type each time on the same fiber  
(which is the proper transform of the line $Q$) to reach $F_{n+1}$.  
Then we apply $n$ elementary transformations of the second type each time 
on the same fiber (which is the proper transform of the line $P$) to return back  
to $F_1$. Then we blow down the exceptional section (whose
self-intersection is now $-1$), and we get again  
$\P^2$.  This process defines a family of Cremona transformations 
from $\P^2$ to $\P^2$.  
 
For each $n$, we get a different Cremona transformation. Uluda\u g has 
shown in \cite{uludag1} that applying to a curve $C$ a  
Cremona transformation whose process reaches $F_{n+1}$,  
yields a new curve $\tilde C$ such that its fundamental group  
$\pi _1 (\P^2-\tilde C)$  
is a central extension of $\pi _1 (\P^2 -C)$ by a cyclic group of  
order $n+1$.

\section{Generalizations of Uluda\u g's method}\label{main} 
In this section we present some generalizations for Uluda\u g's 
method. These generalizations yield new ways 
to construct curves with deep singularities whose fundamental groups 
can be controlled, though they produce no more new groups as 
fundamental groups than the original method of Uluda\u g.  
 
In the first step, we generalize Uluda\u g's method in the following way: 
instead of using the same fiber all the time to perform the 
elementary transformations of the first type, we will use 
an arbitrary finite number of different fibers for performing them 
(Subsection \ref{second_step}).  
In Subsection \ref{general2}, we generalize our construction 
more, and we allow an arbitrary finite number of fibers for
performing the elementary transformations of the second type too. 
In Subsection \ref{special_case_section} we discuss a particularly
interesting special case, where we perform all the transformations (of
both types) on the same fiber.

Before we pass to the generalizations  
of the methods and their proofs, we first have to introduce 
{\it meridians} and prove a lemma which will be needed in the sequel 
(Subsection \ref{meridian}). Also, we have to understand what happens
to the fundamental group when we glue a line back to $\P^2$
(Subsection \ref{fund_group}).  
 
\subsection{Meridians and a generalization of Fujita's lemma}\label{meridian} 
\begin{figure}[ht!]
\epsfysize=4cm 
\centerline{\epsfbox{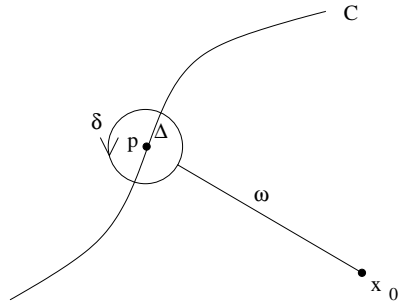}} 
\caption{A meridian of a curve $C$ at a point $p$}\label{meridian_fig} 
\end{figure}

As in Uluda\u g's proof \cite{uludag1}, in order to find the relations 
induced by the additional lines, we have to calculate the 
{\it meridians} of these lines.  
We first recall the definition of a meridian of a curve $C$ at a point
$p$ (see \cite{uludag1,uludag2}):  
Let $\De$ be a smooth analytical branch meeting $C$ transversally 
at $p$ and let $x_0 \in \P^2-C$ be a base point. 
Take a path $\omega$ joining $x_0$ to a boundary point of  
$\De$, and define {\it the meridian of $C$ at $p$} to be the loop  
$\mu _p = \omega \cdot \de \cdot \omega^{-1}$, where $\de$ is 
the boundary of $\De$, oriented in the positive sense 
(see Figure \ref{meridian_fig}).

For computing the meridians in our case, one has to  use some rules. 
The first rule deals with the connection between the meridian 
of the curve $C$ at a point $p \in C$ before the blow up 
and the meridian of the exceptional section created by the blow up 
(see \cite[p. 5]{uludag1}): 
 
\begin{claim}\label{meridian_exceptional} 
Let $\si_p: X\to \P^2$ be a blow up of the point $p \in C$, and 
let $E \subset X$ be the exceptional section. Let $C'=\si_p^{-1}(C)$. 
Then, the loop $\si_p^{-1} (\mu _p)$ is 
the meridian of $E$ at a smooth point $q \in E-C'$. 
\end{claim} 
 
\begin{figure}[ht!]
\epsfysize=5cm 
\centerline{\epsfbox{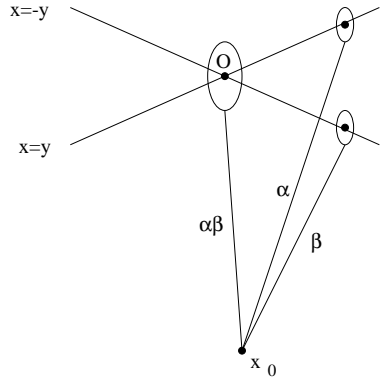}} 
\caption{The situation of Fujita's lemma}\label{fujita_fig} 
\end{figure}

The second rule deals with the meridian at a nodal point: 
\begin{lemma}[Fujita {\cite[p. 540, Lemma 7.17]{fujita}}] \label{fujita}  
Let $B$ be a ball centered at the origin $O$ of $\C^2$,  
and consider the curve $C$ defined by $x^2-y^2=0$. $C$ has 
an ordinary double point at the origin and $\pi_1(B-C)=\Z^2$. 
Take meridians $\alpha$ and $\beta$ of $C$  
on the branches $x=y$ and $x=-y$ at smooth points respectively.  
Then $\alpha \beta$ is a meridian of $C$ at the node $O$ (see Figure 
\ref{fujita_fig}).\end{lemma}

For our generalizations, we need the following more 
general version of this lemma: 

\begin{lemma} \label{gen_fujita} 
Let $B$ be a ball centered at the origin $O$ of $\C^2$,  
and consider the curve $C$ defined by $\prod_{i=1}^k (y-m_i x)=0$
where the $m_i$ are some complex numbers, ordered by their distinct 
arguments, and $k \geq 3$. The curve $C$ has 
an ordinary singular point of multiplicity $k$ at the origin and 
$\pi_1(B-C)=\Z \oplus \F_{k-1}$, where $\F_{k-1}$ is the free group 
on $k-1$ generators. Take meridians $\alpha_i$ of $C$  
on the branches $y=m_i x$ at smooth points respectively.  
Then $\alpha_1 \alpha_2 \cdots \alpha_k$ is a meridian of $C$ 
at the intersection point $O$.    
\end{lemma} 
 
\begin{proof} 
The fact that $\pi_1(B-C)=\Z \oplus \F_{k-1}$ is from \cite{GaTe} or \cite{randell}: Using 
van Kampen's method \cite{vK}, the presentation of the group $\pi_1(B-C)$ is: 
$$\langle \alpha_1, \alpha_2, \cdots, \alpha_k \ | \ \alpha_1 \alpha_2 \cdots \alpha_k=\alpha_2 \cdots \alpha_k \alpha_1 = \cdots = \alpha_k \alpha_1 \cdots \alpha_{k-1} \rangle,$$ 
which can be written also as: 
$$\langle \alpha, \alpha_2, \cdots, \alpha_k \ | \ \alpha \alpha_2 = \alpha_2 \alpha, \alpha \alpha_3 = \alpha_3 \alpha, \cdots ,\alpha \alpha_k = \alpha_k \alpha \rangle$$ 
where $\alpha=\alpha_1 \alpha_2 \cdots \alpha_k$ (see \cite{GaTe}). This presentation  
can be written also as: 
$$\langle \alpha \rangle \oplus \langle \alpha_2, \cdots, \alpha_k \rangle \cong \Z \oplus \F_{k-1},$$ 
and hence the generator of the cyclic group (which is also the center
of the group) is $\alpha=\alpha_1 \alpha_2 \cdots \alpha_k$. 
 
In order to prove the lemma, we have to show that this generator 
is indeed the meridian of the exceptional section $E$ after 
we blow up the point $O$. When we blow the point $O$ up, we get  
Hirzebruch surface $F_1 \to \P^1$.
After deleting the exceptional section and another disjoint section
(which corresponds to the line at infinity), 
we also throw away $k$ fibers corresponding to the $k$ lines which 
passed through $O$ before the blow up. The resulting {\it affine surface} can
be decomposed as a product:  
$(\C -\{ {\rm pt} \}) \times (\P^1- \{ {\rm k\ points} \})$.  Hence, 
its fundamental group can be decomposed into a direct sum too:  
$$\pi_1(\C - \{ {\rm pt} \} ) \oplus \pi_1(\P^1-\{{\rm k\ points}\}) \cong \Z \oplus \F_{k-1}.$$ 
Now, since the cyclic group $\Z$ is in the center of the group $\Z \oplus \F_{k-1}$, its generator 
corresponds indeed to the meridian of the base $E$. Due to the fact
that the generator of the cyclic group in the center 
is $\alpha=\alpha_1 \alpha_2 \cdots \alpha_k$, 
therefore since both the meridian and the generator 
$\alpha=\alpha_1 \alpha_2 \cdots \alpha_k$ generate the infinite 
cyclic group which is the center of the group, this generator is a
meridian of the exceptional section $E$, 
and a meridian of the curve $C$ at the intersection point $O$ too.    
\end{proof}

\subsection{The effect on the fundamental group while gluing back a line}\label{fund_group}

In this short subsection, we prove a simple but useful lemma about the 
effect on the fundamental group when we glue a line back to $\P^2$.

Zaidenberg \cite{Zai} has proved the following lemma:
\begin{lemma}[{Zaidenberg \cite[Lemma 2.3(a)]{Zai}}]
Let $D$ be a closed hypersurface in a complex manifold $M$. Then the
group ${\rm Ker} \{ i_*: \pi_1(M-D) \to \pi_1(M) \}$ is generated by
the vanishing loops of $D$. In particular, if $D$ is irreducible, then,
as a normal subgroup, this group is generated by any of these loops.
\end{lemma}

Let $C$ be a plane curve. 
Substituting $\P^2-C$ for $M$ and a line $L$ for $D$, we get that
${\rm Ker} \{ i_*: \pi_1(\P^2-(C \cup L)) \to \pi_1(\P^2-C) \}$ is
generated by the vanishing loops (=meridians) of $L$. Since $L$ is a
line, we have:
$${\rm Ker} \{ i_*: \pi_1(\P^2-(C \cup L)) \to \pi_1(\P^2-C) \} =
\<\<\mu\>\>$$
where $\mu$ is a meridian of $L$.

Therefore, it is easy to deduce the following lemma:
\begin{lemma}\label{fund_group_lemma}
Let $C$ be a plane projective curve and let $L$ be a line. Let $\mu$
be a meridian of $L$. Then:
$$\pi_1(\P^2-(C \cup L))/ \< \mu \> \cong \pi_1(\P^2-C)$$  
\end{lemma}

\subsection{The first generalization}\label{second_step} 
The first generalization of Uluda\u g's construction is as follows: 
Let $C$ be the initial 
plane curve, and let $n_1, \cdots, n_k$ be $k$ given natural numbers.  
In this construction, we have $k$ different lines $Q_1, \cdots, Q_k$ 
which all meet $C$ transversally. They intersect in a point $O$ 
outside $C$, in such a way that if we put a disk $D$ centered at $O$, 
the intersection points of the lines and the boundary of $D$ are organized 
counterclockwise on the boundary of $D$.  
The additional line $P$ passes via $O$ too and intersects $C$ 
transversally.     
Now we blow up the point $O$, in order to get Hirzebruch surface 
$F_1$. Then we apply $n_i$ elementary transformations of the first type on 
the proper transform of $Q_i$ for all $i=1,\dots, k$. After 
this step, we have reached Hirzebruch surface 
$F_{(\sum_{i=1}^k n_i)+1}$. Now, we apply $\sum_{i=1}^k n_i$  
elementary transformations of the second type on the proper transform of 
$P$. Then, we reach back $F_1$. At last, we blow down the exceptional 
section, and we get again $\P^2$. This defines a family of   
Cremona transformations from $\P^2$ to $\P^2$.  
 
For any $k$-tuple $(n_1,\cdots,n_k) \in \N^k$, we get a different Cremona 
transformation, denoted by $T_{(n_1,\cdots, n_k)}$. We 
will show that the new curve $\tilde C = T_{(n_1,\cdots, n_k)} (C)$ has a  
fundamental group $\pi _1 (\P^2-\tilde C)$  
which is a central extension of $\pi _1 (\P^2 -C)$  
by a finite cyclic group of order $(\sum_{i=1}^k n_i)+1$.

\begin{rem}\label{rem_uludag_successive}
Before formulating the result, we note that
the curve $\tilde C$ obtained by the transformation $T_{(n_1,\cdots, n_k)}$ 
can not be obtained by any successive
applications of Uluda\u g's original method, since for instance two appropriate 
applications will yield an extension of order $n_1+n_2+2$ whereas the extension of the
Cremona transformation $T_{(n_1,n_2)}$ (for $k=2$) is of order $n_1+n_2+1$. Also the
obtained singularities will be different 
(see Section \ref{extension_cyclic}).
\end{rem}

\begin{thm} \label{gen_2step} 
Let $C$ be a plane projective curve and $G=\pi _1 (\P ^2 -C)$. 
Then for any $k$-tuple $(n_1,\cdots,n_k)\in \N^k$, 
the curve $\tilde C = T_{(n_1,\cdots, n_k)} (C)$ is 
birational to $C$, and its fundamental group 
$\tilde G=\pi _1 (\P ^2 -\tilde C)$ is a  
central extension of $G$ by $\Z/((\sum_{i=1}^k n_i)+1)\Z$: 
$$1 \rightarrow \Z/((\sum_{i=1}^k n_i)+1)\Z \rightarrow \tilde G \rightarrow G \rightarrow 1$$ 
Moreover, if $C$ has $r$ irreducible components so is $\tilde C$.   
\end{thm}  
 
\begin{proof} 
We start with the observation that  since  
the blow up and the blow down transformations are birational transformations and the fact
that the Cremona transformation is composed of blow up and blow down transformations,  
we have that $\tilde C$ is birational to $C$, and  
the number of irreducible components is preserved by the construction.  
 
Let $\tilde P=T_{(n_1,\cdots, n_k)}(P)$ and 
$\tilde Q_i= T_{(n_1,\cdots, n_k)}(Q_i)$ for $1 \leq i \leq k$
(one should notice here that $P$ is not actually mapped to $\tilde P$, since 
the transformation blows it down. We choose the new fiber that replaces the
original one).
This Cremona transformation defines an isomorphism: 
$$\P^2 - (C \cup P \cup (\bigcup_{i=1}^k Q_i)) \cong \P^2 - (\tilde C \cup \tilde P \cup (\bigcup_{i=1}^k \tilde Q_i)),$$ 
which induces an isomorphism of the corresponding fundamental groups: 
$$\pi_1(\P^2 - (C \cup P \cup (\bigcup_{i=1}^k Q_i))) \cong \pi_1(\P^2 - (\tilde C \cup \tilde P \cup (\bigcup_{i=1}^k \tilde Q_i))).$$ 
In order to compute $\pi_1(\P^2 - \tilde C)$, we have to add 
the relations correspond 
to gluing back the lines  $\tilde P$ and $\tilde Q_i, 1 \leq i \leq k$. 
 
Let $\beta$ and $\alpha_i$ be the meridians of the lines $P$ and $Q_i$ at smooth points
respectively. Using Lemma \ref{gen_fujita}, we get that the meridian of the  
curve $C \cup P \cup (\bigcup_{i=1}^k Q_i)$ at the point $O$ 
(which is the intersection of $k+1$ lines: 
$P,Q_1, \cdots, Q_k$) is  
$\beta \alpha_1 \cdots \alpha_k$, and hence the meridian of $E$,  
which is the blow up of this point, is $\beta \alpha_1 \cdots
\alpha_k$ too.  
    
We now compute the meridian of $\tilde Q_i$ ($1 \leq i \leq k$).
In each elementary transformation of the first type on $Q_i$ and 
its proper transforms, 
we blow up the intersection point between the fiber and the exceptional section, 
and then we blow down the fiber.  Hence, by Fujita's lemma (Lemma \ref{fujita}), 
each elementary transformation multiplies the current meridian 
of the proper transform of $Q_i$ by the meridian 
of the exceptional section.
Therefore, we finally get that the meridian of $\tilde Q_i$ is
$(\beta \alpha_1 \cdots \alpha_k)^{n_i} \cdot \alpha_i$
(see Figure \ref{meridians1} for the effect of one elementary transformation 
of the first type where $k=2$). 
The elementary transformations on $Q_i$ and its proper transforms do not change 
the meridians of the other lines.

\begin{figure}[ht!] 
\epsfysize=8cm 
\centerline{\epsfbox{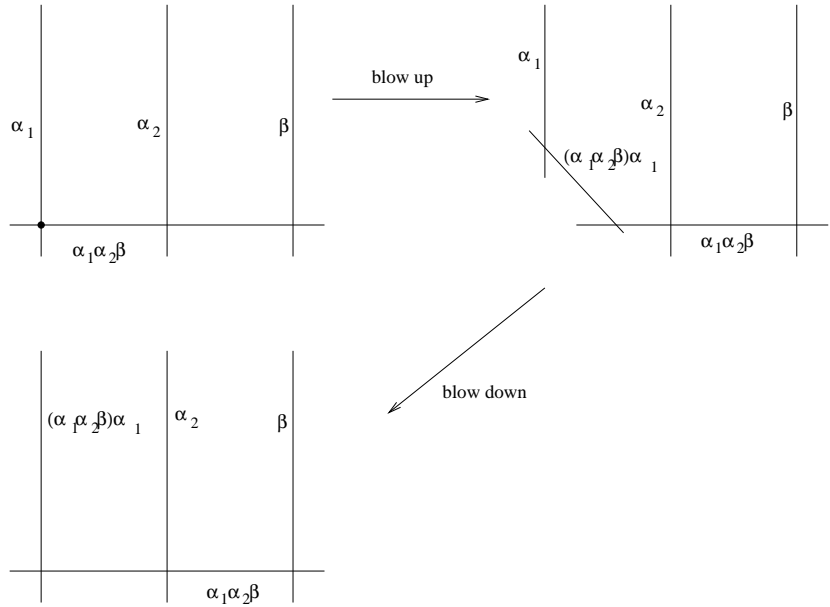}} 
\caption{The effect of an elementary transformation of the first type
  on the meridians} \label{meridians1}
\end{figure}                                         

We have to compute the meridian of $\tilde P$ too. 
In each elementary transformation of the second type on $P$ and 
its proper transforms, 
we blow up a point on the fiber, and then we blow down the fiber.  
Hence, this type of elementary transformations preserves the current meridian 
and we get that the meridian of $\tilde P$ remains $\beta$
(see Figure \ref{meridians2} for the effect of one elementary transformation 
of the second type where $k=2$). 
The elementary transformations on $P$ and its proper transforms do not change 
the meridians of the other lines.

\begin{figure}[ht!] 
\epsfysize=8cm 
\centerline{\epsfbox{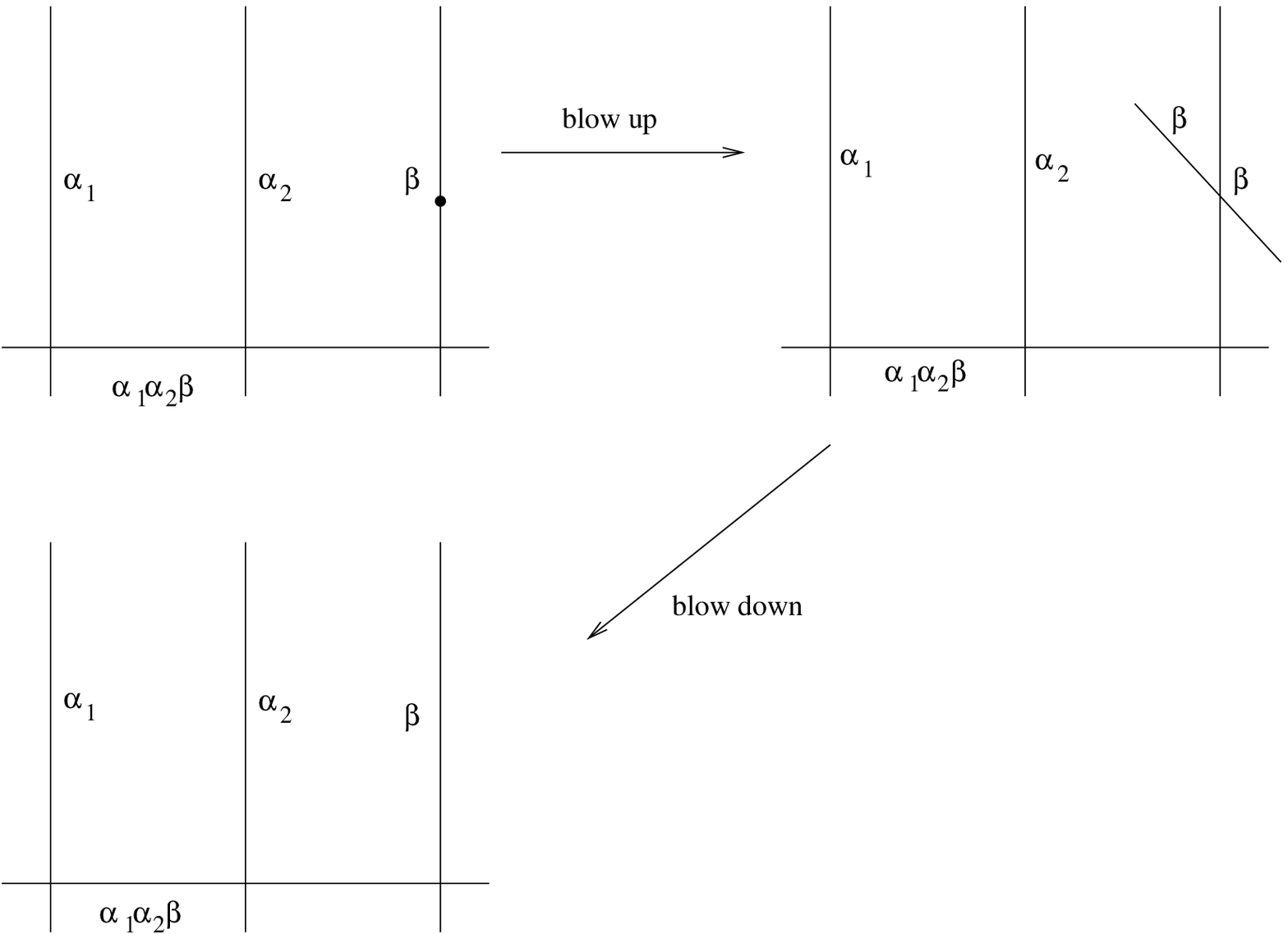}} 
\caption{The effect of an elementary transformation of the second type
  on the meridians} \label{meridians2}
\end{figure}                                         

Therefore, by Lemma \ref{fund_group_lemma}, we conclude that:

$\pi_1 (\P^2- \tilde C) \cong$\vspace{-10pt}
$$\pi_1(\P^2-(\tilde C \cup \tilde P \cup  (\bigcup_{i=1}^k \tilde Q_i))) / \langle \beta,(\beta \alpha_1 \cdots \alpha_k)^{n_1} \alpha_1, \cdots, (\beta \alpha_1 \cdots \alpha_k)^{n_k} \alpha_k \rangle,$$ 
which is equivalent to: 
$$\pi_1 (\P^2- \tilde C) \cong \pi_1(\P^2-(\tilde C \cup \tilde P \cup  (\bigcup_{i=1}^k \tilde Q_i))) / \langle \beta,(\alpha_1 \cdots \alpha_k)^{n_1} \alpha_1, \cdots, (\alpha_1 \cdots \alpha_k)^{n_k} \alpha_k \rangle$$ 
since the connecting relations between $\beta,\alpha_1, \cdots, \alpha_k$ are  
$$\beta \alpha_1 \cdots \alpha_k =\alpha_1 \cdots \alpha_k \beta = \cdots = \alpha_k \beta \alpha_1 \cdots \alpha_{k-1}$$  
(see the proof of Lemma \ref{gen_fujita}), because all the $k+1$ lines 
intersect in one point $O$. 
 
Since
$$\P^2 - ( C \cup P \cup  (\bigcup_{i=1}^k Q_i)) \cong \P^2 - (\tilde C \cup \tilde P \cup (\bigcup_{i=1}^k \tilde Q_i)),$$ 
we can also write: 

$$\pi_1 (\P^2- \tilde C) \cong \pi_1(\P^2-(C \cup P \cup (\bigcup_{i=1}^k Q_i))) / \langle \beta,(\alpha_1 \cdots \alpha_k)^{n_1} \alpha_1, \cdots, (\alpha_1 \cdots \alpha_k)^{n_k} \alpha_k \rangle$$ 
When we pass to the quotient by $\langle \beta \rangle$ which 
corresponds to gluing the line  
$P$ back into $\P^2$ in the original configuration 
(by Lemma \ref{fund_group_lemma}), we get that: 
$$\pi_1 (\P^2- \tilde C) \cong \pi_1(\P^2-(C \cup (\bigcup_{i=1}^k Q_i))) / \langle (\alpha_1 \cdots \alpha_k)^{n_1} \alpha_1, \cdots, (\alpha_1 \cdots \alpha_k)^{n_k} \alpha_k \rangle.$$ 
While moving to this quotient, we get the following 
$k$ cyclic relations for the subgroup in the denominator:  
$$\alpha_1 \alpha_2 \cdots \alpha_k = \alpha_2 \cdots \alpha_k \alpha_1 = \cdots = \alpha_k \alpha_1 \alpha_2 \cdots \alpha_{k-1},$$  
since $\beta=1$ in the quotient. These relations can 
be presented also as the following set of relations, where  
$\alpha =\alpha_1 \alpha_2 \cdots \alpha_k$: 
$$\{ \alpha \alpha_2 =\alpha_2 \alpha, \alpha \alpha_3 =\alpha_3 \alpha, \cdots , \alpha \alpha_k =\alpha_k \alpha \}$$  
By Lemma \ref{fund_group_lemma}, when we take the quotient of 
$\pi_1(\P^2-(C \cup (\bigcup_{i=1}^k Q_i)))$ by the normal subgroup generated
by the meridians of $Q_1, \cdots, Q_k$, we get $\pi_1(\P^2-C)$. 
Due to the
fact that all the lines are intersected at $O$, then the subgroup
generated by the meridians $Q_1, \cdots, Q_k$ is of the form 
(see the proof of Lemma \ref{gen_fujita}):
$$\Z \oplus \F_{k-1} \cong \langle \alpha, \alpha_2, \cdots, \alpha_k \ | \ \alpha\alpha_i=\alpha_i\alpha, 2 \leq i \leq k \rangle$$ 
where $\alpha = \alpha_1 \cdots \alpha_k$. 

We have
$${{\pi_1(\P^2-(C\cup (\bigcup_{i=1}^k Q_i)))/ \langle \alpha^{n_1} \alpha_1, \cdots, \alpha^{n_k} \alpha_k \rangle} \over {\langle \alpha, \alpha_2, \cdots, \alpha_k \ | \ \alpha\alpha_i=\alpha_i\alpha, 2 \leq i \leq k \rangle / \langle \alpha^{n_1} \alpha_1, \cdots, \alpha^{n_k} \alpha_k \rangle }} \cong $$ 
$$\cong \pi_1(\P^2-(C\cup (\bigcup_{i=1}^k Q_i))) / \langle \alpha, \alpha_2, \cdots, \alpha_k \rangle$$ 
where $\alpha = \alpha_1 \cdots \alpha_k$, and therefore: 
$${{\pi_1(\P^2-\tilde C) } \over {\left( {\langle \alpha, \alpha_2, \cdots, \alpha_k \ | \ \alpha\alpha_i=\alpha_i\alpha, 2 \leq i \leq k \rangle \over \langle \alpha^{n_1} \alpha_1, \cdots, \alpha^{n_k} \alpha_k \rangle} \right) }} \cong \pi_1(\P^2-C)$$ 
which can be written as the following extension: 
$$1 \to {\langle \alpha, \alpha_2, \cdots, \alpha_k \ | \ \alpha\alpha_i=\alpha_i\alpha, 2 \leq i \leq k \rangle \over \langle \alpha^{n_1} \alpha_1, \cdots, \alpha^{n_k} \alpha_k \rangle } \to \pi_1(\P^2-\tilde C) \to \pi_1(\P^2-C) \to 1$$
Hence, it remains to show that this extension is central and that 
$${\langle \alpha, \alpha_2, \cdots, \alpha_k \ | \ \alpha\alpha_i=\alpha_i\alpha, 2 \leq i \leq k \rangle \over \langle \alpha^{n_1} \alpha_1, \cdots, \alpha^{n_k} \alpha_k \rangle} \cong \Z / ((\sum_{i=1}^k n_i)+1)\Z$$ 
The centrality of the extension is a little bit 
tricky. Although $\alpha_i$ does not commute with $\alpha_j$ 
($2 \leq i,j \leq k$, $i \not = j$) in $\pi_1(P^2-(C\cup (\bigcup_{i=1}^k Q_i)))$, the generators of the quotient group 
$\tilde G =\pi_1(\P^2-\tilde C)$ are the meridians 
of $C$ and only one more generator - $\alpha$ (the other generators - 
$\alpha_2, \cdots, \alpha_k$ - corresponding to the lines 
$Q_2,\cdots, Q_k$ in the bigger group 
$\pi_1(\P^2-(C\cup (\bigcup_{i=1}^k Q_i)))$ disappear in $\tilde G$ by the 
additional relations $\alpha^{n_2} \alpha_2, \cdots, \alpha^{n_k} \alpha_k$). 
This generator indeed commutes 
with the meridians of $C$ in $\pi_1(\P^2-C)$ (as it is equal to 
the multiplication of all the $\alpha_i$'s, which commute 
with the meridians of $C$ in $\pi_1(\P^2-C)$ since we choose the lines
$Q_1,\cdots, Q_k$ to be all transversally intersected with $C$), 
and hence we get that the extension is central.

Now, we show that: 
$${\langle \alpha, \alpha_2, \cdots, \alpha_k \ | \ \alpha\alpha_i=\alpha_i\alpha, 2 \leq i \leq k \rangle \over \langle \alpha^{n_1} \alpha_1, \cdots, \alpha^{n_k} \alpha_k \rangle} \cong \Z / ((\sum_{i=1}^k n_i)+1)\Z$$ 
Since  
$$\alpha^{n_1} \alpha_1 \cdots \alpha^{n_k} \alpha_k = \alpha^{n_1+ \cdots +n_k}(\alpha_1 \cdots \alpha_k) = \alpha^{n_1+ \cdots +n_k} \alpha =\alpha^{n_1+ \cdots +n_k+1},$$  
we can change the presentation of the subgroup in the denominator to the following: 
$$\langle \alpha^{n_1+ \cdots +n_k+1},\alpha^{n_2} \alpha_2, \cdots, \alpha^{n_k} \alpha_k \ | \  \alpha\alpha_i=\alpha_i\alpha, 2 \leq i \leq k \rangle.$$ 
On the other hand, one can see that the following is 
a presentation of $\Z \oplus \F_{k-1}$, which is the numerator: 
$$\Z \oplus \F_{k-1} \cong  \langle  \alpha ,\alpha^{n_2} \alpha_2, \cdots, \alpha^{n_k} \alpha_k \ | \  \alpha\alpha_i=\alpha_i\alpha, 2 \leq i \leq k \rangle,$$ 
since all the $\alpha_i, 2 \leq i \leq k,$ can be achieved by 
the new set of generators.  
By these two new presentations, one can easily see that 
the quotient of these two groups is $\Z / ((\sum_{i=1}^k n_i)+1)\Z$ as needed, 
and hence Theorem \ref{gen_2step} is proved. 
\end{proof} 

\subsection{A slightly more general construction}\label{general2}
During the proof of the last theorem (before Figure \ref{meridians2}), 
we have shown that elementary transformations of
the second type do not affect the meridians of the fibers which we perform
the transformations on. Therefore, we can generalize our
construction to the following one: instead of performing all the
elementary transformations of the second type on the same fiber $P$, we can
apply them on several fibers $P_1,\cdots, P_t$, with the condition that the
total number of applications of elementary transformations 
of the second type will be equal to the total number of applications of 
elementary transformations of the first type.    

Using this observation, we can describe a slightly more general construction:
Let $C$ be the initial plane curve, and let $n_1, \cdots, n_k$ and 
$m_1, \cdots, m_l$ be two sets of $k$ and $l$ given natural numbers, such that 
$\sum n_i = \sum m_j$.  We start with $k+l$ different lines 
$Q_1,\cdots, Q_k$ and $P_1,\cdots, P_l$  which all meet $C$ transversally. 
They all intersect in a point $O$ outside $C$, in such a way 
that if we put a disk $D$ centered at $O$, the intersection points of 
the lines and the boundary of $D$ are organized counterclockwise on the boundary 
of $D$.
Now we blow up the point $O$, in order to get Hirzebruch surface 
$F_1$. Then we apply $n_i$ elementary transformations of the first type on 
the proper transform of $Q_i$ for all $i=1,\dots, k$. After 
this step, we have reached Hirzebruch surface 
$F_{(\sum_{i=1}^k n_i)+1}$. Now, we apply $m_j$  
elementary transformations of the second type on the proper transform of 
$P_j$ for all $j=1,\dots, l$. Since $\sum n_i = \sum m_j$, 
we reach back $F_1$. At last, we blow down the exceptional 
section (whose self-intersection is now $-1$), 
and we get again $\P^2$. This defines a family of   
Cremona transformations from $\P^2$ to $\P^2$.  
 
For any $(k+l)$-tuple $(n_1,\cdots,n_k,m_1,\cdots,m_l) \in \N^{k+l}$ such that 
$\sum n_i = \sum m_i$, we get a Cremona transformation, 
denoted by $T_{(n_1,\cdots, n_k;m_1,\cdots,m_l)}$. 
Then we can state:
  
\begin{prop} \label{gen_2step_coro} 
Let $C$ be a plane projective curve and $G=\pi _1 (\P ^2 -C)$. 
Then for any $(k+l)$-tuple $(n_1,\cdots,n_k,m_1,\cdots,m_l)\in \N^{k+l}$
such that $\sum n_i = \sum m_i$, the curve 
$\tilde C = T_{(n_1,\cdots, n_k;m_1,\cdots,m_l)} (C)$  
is birational to $C$ and its fundamental group 
$\tilde G=\pi _1 (\P ^2 -\tilde C)$ is a  
central extension of $G$ by $\Z/((\sum_{i=1}^k n_i)+1)\Z$: 
$$1 \rightarrow \Z/((\sum_{i=1}^k n_i)+1)\Z \rightarrow \tilde G \rightarrow G \rightarrow 1$$ 
Moreover, if $C$ has $r$ irreducible components so is $\tilde C$.   
\end{prop}  

\subsection{An interesting special case}\label{special_case_section}

Just before finishing the section of the
constructions, we want to concentrate on an interesting special
construction. In this construction, we perform all the elementary
transformations from both types on the same fiber. 

First, we define this construction precisely, and then we prove that
the fundamental group of the resulting curve is again a central extension of
the original curve, as we had in the previous constructions. We have
to prove it, since the proof is slightly different from the proof of
the previous case.

We start with a curve $C$ and one additional line $L$ in $\P^2$ which
intersects $C$ transversally.
Blow up a point $O$ on $L$ (which does not belong to $C$) in order to reach
$F_1$. Then perform $n$ elementary transformations of the first type
on the proper transform of $L$. Hence we reach $F_{n+1}$. Now, perform $n$
elementary transformations of the second type again on the proper
transform of $L$. Now, blow down the exceptional section (which now has 
self-intersection $-1$) in order to return to $\P^2$. This construction
defines a Cremona transformation from $\P^2$ to $\P^2$, which we
denote by $T_n$. Then we have the following result:

\begin{prop} \label{special_case} 
Let $C$ be a plane projective curve and $G=\pi _1 (\P ^2 -C)$. 
Then for any natural number $n\in \N$, the curve 
$\tilde C = T_n (C)$ is birational to $C$ and its fundamental group 
$\tilde G=\pi _1 (\P ^2 -\tilde C)$ is a  
central extension of $G$ by $\Z/(n+1)\Z$: 
$$1 \rightarrow \Z/(n+1)\Z \rightarrow \tilde G \rightarrow G \rightarrow 1$$ 
Moreover, if $C$ has $r$ irreducible components so is $\tilde C$.   
\end{prop}  

\begin{proof} 
As before, $\tilde C$ is birational to $C$, and number of irreducible components
is preserved by $T_n$.
 
Let $\tilde L = T_n(L)$. $T_n$ induces an isomorphism 
of the corresponding fundamental groups: 
$$\pi_1(\P^2 - (C \cup L)) \cong \pi_1(\P^2 - (\tilde C \cup \tilde L))$$ 
Let $\alpha$ be the meridian of $L$ at a smooth point.
By Claim \ref{meridian_exceptional}, we get that the meridian of the
exceptional section $E$, which is the blow up of the point $O \in L$, 
is again $\alpha$.
    
After we apply the sequence of $n$ elementary transformations 
of the first type using the fiber $L$ and its proper transforms, we have by 
Fujita's lemma (Lemma \ref{fujita}) that the meridian of the proper transform of $L$
is $\alpha ^{n+1}$.
As before, the applications of elementary 
transformations of the second type  
and the final blow down do not change this meridian. 
 
Therefore, by Lemma \ref{fund_group_lemma}, we conclude that: 
$$\pi_1 (\P^2- \tilde C) \cong \pi_1(\P^2-(\tilde C \cup \tilde L)) / \langle \alpha^{n+1} \rangle$$ 
As before, it is equivalent to:
$$\pi_1 (\P^2- \tilde C) \cong \pi_1(\P^2-(C \cup L)) / \langle \alpha^{n+1} \rangle.$$
By Lemma \ref{fund_group_lemma}, 
$\pi_1(\P^2-(C\cup L))/\langle \alpha \rangle =\pi_1(\P^2-C)$. 
Hence:
$$\pi_1(\P^2-\tilde C)/ \langle \alpha \rangle \cong \pi_1(\P^2-C),$$
which can be written as the following extension: 
$$1 \to (\langle \alpha \rangle/\langle
\alpha^{n+1} \rangle) \to \pi_1(\P^2-\tilde C) \to \pi_1(\P^2-C) \to 1.$$ 
Obviously, $\langle \alpha \rangle/\langle\alpha^{n+1} \rangle \cong
\Z/(n+1)\Z$, and since $\alpha$ commutes with the generators 
of $\pi_1(\P^2-C)$ (since the line $L$ intersects $C$ transversally), 
the extension is central.  
\end{proof}

\section{Properties of groups preserved by the constructions}\label{properties}  
In this section, we indicate some properties of the fundamental group 
which are preserved by the constructions of the previous section.

We start with an interesting property about the splitness of the
central extension we have in the constructions. Using this property,   
we will show the following important property of the fundamental group
of the resulting curve (Proposition \ref{cyclicity}): 
if we start with an irreducible 
curve which has a cyclic group as the fundamental group 
(such as smooth irreducible curves), then the resulting fundamental group
will be cyclic too. The importance of this property is that although the
constructions add to the curve deep singularities (as is proved in 
Proposition \ref{sing_general}), the fundamental group of the curve is
still cyclic. Hence, these constructions may yield families of plane curves
which have some deep singularities but have cyclic fundamental groups.  

\begin{prop}\label{cyclicity_general} 
Let $C$ be a plane curve with $r$ irreducible components.
Let $n \in \N$. Let $\tilde C$ be the curve whose 
fundamental group $\tilde G$ is obtained 
from  $G=\pi _1 (\P ^2 -C)$ by a central extension by $\Z/(n+1)\Z$.
If the abelian group $H_1 (\P ^2 -C)$ has $r$ direct summands and the orders 
of the summands are not coprime to $n+1$, then the extension does not split, 
i.e. $\tilde G \not \cong G \oplus \Z/(n+1)\Z$.
\end{prop} 

\begin{rem}
In case of coprimeness, we obviously have 
$$\Z/(mn)\Z \cong \Z/m\Z \oplus \Z/n\Z,$$
and this is the reason for ruling this case out in the proposition.
\end{rem}

\begin{proof}[Proof of Proposition \ref{cyclicity_general}] 
On the contrary, assume that $\tilde G \cong G \oplus \Z/(n+1)\Z$.   
As $H_1(X)$ is the abelinization of $\pi _1 (X)$,
we have that 
$$H_1 (\P ^2 -\tilde C) \cong {\rm Ab}(\tilde G) \cong {\rm Ab}(G \oplus
\Z/(n+1)\Z) \cong H_1 (\P ^2 -C) \oplus \Z/(n+1)\Z.$$
Since $n+1$ is not coprime to the orders of the summands of $H_1 (\P ^2 -\tilde C)$, 
$H_1 (\P ^2 -\tilde C)$ has $r+1$ proper direct summands.
This contradicts the fact that $\tilde C$ has only $r$ 
irreducible components, as the number of irreducible components 
is preserved by the constructions.  
\end{proof}

\begin{prop}\label{cyclicity} 
Let $C$ be an irreducible plane curve with 
a cyclic fundamental group $\Z/r\Z$.  
Let $n \in \N$. Let $\tilde C$ be the curve whose fundamental 
group $\tilde G$ is obtained 
from  $G=\pi _1 (\P ^2 -C)$ by a central extension by $\Z/(n+1)\Z$.  
Then $\tilde G=\pi _1 (\P ^2 -\tilde C)$ is also cyclic of order $r(n+1)$.   
\end{prop} 
 
\begin{proof} 
As $H_1(X)$ is the abelinization of $\pi _1 (X)$, 
we have that $G = H_1 (\P ^2 - C)$. Since $G$ is cyclic,
$H_1 (\P ^2 - C)$ has one direct summand, which equals the number of
irreducible components in $C$ (one too). 
By the previous proposition, the extension does not split, and 
$H_1 (\P ^2 - \tilde C)$ is cyclic too.
 
Since the extension is central, we have that $\tilde G$ 
is abelian of order $r(n+1)$. Hence, $\tilde G=H_1 (\P ^2 - \tilde C)$.
Therefore, $\tilde G$ is cyclic of order $r(n+1)$. 
\end{proof} 
 
\begin{rem}
The condition that the curve is irreducible is essential, since if we
take a reducible curve with a cyclic fundamental group, it is not
guaranteed that the resulting curve will have a cyclic fundamental
group. For example, if we start with a curve $C$ consists of two
intersecting lines whose fundamental group $\pi_1(\P ^2 - C)=\Z$ is 
cyclic, and we apply on it Uluda\u g's method for $n=1$, we get that the
resulting curve $\tilde C$ has a fundamental group 
$\pi_1(\P ^2 - \tilde C)=\Z \oplus (\Z/2\Z)$ 
(see after Proposition \ref{general_pos}) which is not cyclic. 
\end{rem} 

Let $p$ be a prime number. Then:    
\begin{rem} 
Let $C$ be a plane curve with a fundamental group $G$ which is a $p$-group.  
Let $n \in \N$. Let $\tilde C$ be the curve whose fundamental group 
$\tilde G$ is obtained from  $G=\pi _1 (\P ^2 -C)$ by a central extension 
by $\Z/(n+1)\Z$.  Then: if $n+1=p^l$ for some $l$, 
then $\tilde G=\pi _1 (\P ^2 -\tilde C)$ is also a $p$-group.   
\end{rem} 
 
\begin{proof} 
From the extension, we get that $\tilde G/(\Z/ (n+1)\Z) \cong G$. 
Since $n+1=p^l$, $\Z/(n+1)\Z$ is a $p$-group, and  
since $G$ is also a $p$-group, then $\tilde G$ is a $p$-group too.  
\end{proof} 

\begin{rem} \label{rem_decomposition}
Let $C$ be a plane curve with a finite fundamental group $G$.  
Let $n \in \N$. Let $\tilde C$ be the curve whose fundamental group 
$\tilde G$ is obtained 
from  $G=\pi _1 (\P ^2 -C)$ by a central extension by $\Z/(n+1)\Z$.  
If $(|G|,n+1)=1$, then the fundamental group $\tilde G$ of the
resulting curve is a direct sum of $G$ and $\Z/(n+1)\Z$:
$$\tilde G \cong G \oplus \Z/(n+1)\Z$$
\end{rem} 
 
\begin{proof}
Use Theorem 7.77 of \cite{rotman} that ``if $Q$ is a finite group, $K$
is a finite abelian group, and $(|K|,|Q|)=1$, then an extension $G$ of
$Q$ by $K$ is a semidirect product of $K$ and $Q$'', and since the extension is central, semidirect products become
direct products.  
\end{proof} 
 
In the following proposition, we list some more properties of the
fundamental group which are preserved by the constructions. Before
stating it, we remind some definitions.
A group is called {\it polycyclic} if it has a subnormal
series with cyclic factors. A group is called {\it supersolvable} 
if it has a normal series with cyclic factors. We say that a group $G$ is
{\it nilpotent} if its lower central series reaches $1$ 
(see for example \cite{rotman}). 

\begin{prop}\label{proper}
Let $C$ be a plane curve with a fundamental group $G$.  
Let $n \in \N$. Let $\tilde C$ be the curve whose fundamental group 
$\tilde G$ is obtained 
from  $G=\pi _1 (\P ^2 -C)$ by a central extension by $\Z/(n+1)\Z$.  
Then if $G$ has one of the following properties, the fundamental group
$\tilde G$ of the resulting curve has that property too:
\begin{enumerate}
\item Finite.
\item Non-abelian.
\item Solvable.
\item Supersolvable.
\item Polycyclic.
\item Nilpotent.
\end{enumerate}
\end{prop} 
 
\begin{proof} 
(1-2)\qua Trivial.

(3)\qua From the extension, we get that $\tilde G/(\Z/ (n+1)\Z) \cong G$. 
But it is known \cite[Theorem 5.17]{rotman} that if $\Z/ (n+1)\Z$ and 
$G$ are both solvable, then $\tilde G$ is solvable too. 

(4)\qua From the extension, we have that $\tilde G/(\Z/ (n+1)\Z) \cong G$.
It is easy to show (very similar to the solvable case, 
see \cite[Theorem 5.17]{rotman}) that if $\Z/ (n+1)\Z$ and 
$G$ are both supersolvable, then $\tilde G$ is supersolvable too.

(5)\qua Same proof as (4).

(6)\qua Since $\tilde G$ is a central extension of $G$ by $\Z/ (n+1)\Z$ for a
given $n$, we get that $\tilde G/(\Z/ (n+1)\Z) \cong G$ where   
$\Z/ (n+1)\Z \leq Z(\tilde G)$. 
But it is easy to see that if $\Z/ (n+1)\Z \leq Z(\tilde G)$ 
and $G$ is nilpotent, then $\tilde G$ is nilpotent too 
(see for example \cite[p. 117, Exercise 5.38]{rotman}).
\end{proof} 
 
We note here that if $G$ is a nilpotent group of class $c$ 
(means that the last non-zero term of  
the lower central series is the $c$-th term), then $\tilde G$ is 
a nilpotent group of class $c$ or $c+1$. 
It will be of class $c$ if and only if the 2-cocycle 
defining the extension is symmetric (see \cite{schafer}).

Here we indicate one more family of group properties which are
preserved by the constructions. 

\begin{rem} 
Let $C$ be a plane curve with a fundamental group $G$ which has a
subgroup $N$ of finite index with a special property 
(for example: solvable, nilpotent, etc.).  
Then the fundamental group $\tilde G$ of the resulting curve $\tilde C$ 
has a subgroup $\tilde N$ of finite index with the same special property too. 

In particular, if $G$ is virtually-nilpotent
or virtually-solvable, then $\tilde G$ has the same property as well.  
\end{rem} 

\begin{proof}
Let $n\in \N$. Since $\tilde G$ is an extension of $G$ by 
$\Z/ (n+1)\Z$, then $\tilde G/(\Z/ (n+1)\Z) \cong G$. 
Hence there exists 
$\tilde N \leq \tilde G$ such that $\tilde N/(\Z/ (n+1)\Z) \cong N$.
Obviously $[\tilde G : \tilde N] = [G:N] < \infty$.
Since $\tilde G$ is a central extension, one can show that also  
$\tilde N$ is a central extension of $N$ by $\Z/ (n+1)\Z$ 
(since $\Z/(n+1)\Z$ is a subgroup of the intersection of 
$Z(\tilde G)$ and $\tilde N$, and hence in $Z(\tilde N)$).
Therefore, one can apply Proposition \ref{proper} to
show that the properties of $N$ are moved to $\tilde N$.
\end{proof}

\begin{rem}
Note that all the results of this section hold also for the
constructions of Oka \cite{oka} and Shimada \cite{shimada1}, 
since also in their constructions
the fundamental group of the resulting curve is a central extension
of the fundamental group of the original curve by a cyclic group. 
\end{rem}

\section{The curves obtained by the constructions}\label{extension_cyclic} 
 
In this section we investigate the curves which can  
be obtained using Uluda\u g's original construction and the general constructions
we have presented in the previous sections. 
In the first subsection, we will describe the types 
of singularities which are added  
by the constructions. In Subsection \ref{sec_deg} we compute the
degrees of the resulting curves.  
In the next subsections, we describe some families of curves which  
can be obtained by applying the constructions on several 
different types of curves.

\subsection{The types of singularities which are added to the curves}  
 
At the beginning of this subsection, we want to fix a notation for 
singular points.  
We follow the notations of Flenner and Zaidenberg \cite{zaidenberg}. 
Although in general it is possible that a singular point splits into 
several points by a blow up, in most of our cases (except for 
Proposition \ref{gen_2step_coro}), only one singular point appears by
a blow up.
Hence, any singular point $P$ has a very simple resolution by a sequence of $s$  
blow ups. We denote by $t_i > 1$  $(1 \leq i \leq s)$ 
the multiplicity of the curve at $P$ before the $i$-th 
blow up. Then $[t_1,\cdots, t_s]$ is called {\it the type of the singularity}. 
If we have a sequence of $r$ equal multiplicities $l$, 
we abbreviate it by $l_r$. 
For example, $[2,2,2]=[2_3]$ corresponds to a ramphoid cusp, 
which has to be blown up  
three times (each time of multiplicity $2$) for smoothing it. 

To simplify the description, we also introduce the notion 
of a {\it $d$-tacnode}. 
\begin{defn} 
A {\it $d$-tacnode} is a singular point where $d$ smooth branches of the curve  
have a common tangent at the same point. 
A {\it $d$-tacnode of order $k$} is a $d$-tacnode such that the higher derivatives
of the branches are also equal up to order $k$.
\end{defn} 
  
For example, the usual tacnode is $2$-tacnode of order $1$.

Here we describe the singularities which are added to the curve 
during these constructions:   
 
\begin{prop}\label{sing_general} 
Let $C$ be a curve of degree $d$. Let $(n_1, \cdots, n_k)\in \N^k$ 
be a $k$-tuple.  
Let $\tilde C = T_{(n_1,\cdots, n_k)} (C)$ (see Theorem \ref{gen_2step}). 
 
Then $\tilde C$ has $k+1$ additional singularities to those of $C$:  
$k$ $d$-tacnodes of order $n_i-1$ ($1 \leq i \leq k$), 
and another singular point which  
is a blow down of a $d$-tacnode of order $n_1+\cdots +n_k-1$ 
(i.e., the curve has the following 
additional singularities: $[d_{n_1}], \cdots, [d_{n_k}]$  and 
$[d(n_1+\cdots +n_k),d_{(n_1+\cdots +n_k)}]$). 
\end{prop} 
 
\begin{proof} 
The first blow-up (from $\P^2$ to $F_1$) does not change the curve. 
Each sequence of $n_i$ elementary transformations of the first type 
on the fiber $Q_i$ and  
its proper transforms creates one $d$-tacnode of order $n_i-1$:
The first elementary transformation of the first type 
creates one intersection point (of $d$ branches). 
The second elementary transformation  
of the first type converts it into a $d$-tacnode of order $1$, and  
another $n_i-2$ elementary transformations of the first type convert it 
into a $d$-tacnode of order $n_i-1$.  
This $d$-tacnode is located at the proper transform of $Q_i$, 
outside the exceptional section.  

The sequence of $n_1+\cdots +n_k$ elementary transformations of the 
second type  
(from $F_{n_1+\cdots +n_k+1}$ to $F_1$) creates another 
$d$-tacnode of order $n_1+\cdots +n_k-1$, as in the case of the elementary 
transformations of the first type. The difference is that now this $d$-tacnode is 
located at the exceptional section. Hence, when we blow down this section in order to 
return to $\P^2$, this $d$-tacnode of order $n_1+\cdots +n_k-1$ 
is blown down too to a more complicated singular point. 
\end{proof} 

The following remark states the situation after applying only 
Uluda\u g's original construction.

\begin{rem}\label{sing_uludag} 
Let $C$ be a curve of degree $d$. Let $n \in \N$. Let $\tilde C$ be the curve  
obtained by Uluda\u g's 
original construction for this $n$ (see Section \ref{uludag}).
 
Then $\tilde C$ has two additional singularities to those of $C$: 
a $d$-tacnode of order $n-1$,  
and another singular point which is a blow down  
of a $d$-tacnode of order $n-1$ (i.e., the curve has the following 
additional singularities: $[d_n]$ and $[dn,d_n]$). 
\end{rem}

Here we describe the situation concerning the special case
(Proposition \ref{special_case}).

\begin{prop}\label{sing_special_case} 
Let $C$ be a curve of degree $d$. Let $n \in \N$.
Let $\tilde C = T_n (C)$. 
 
Then $\tilde C$ has one additional singularity to those of $C$:  
a blow down of a $d$-tacnode of order $2n-1$ 
(i.e., the curve has the following 
additional singularity:
$[2nd,d_{2n}]$). 
\end{prop} 
 
\begin{proof} 
As before, the first blow-up does not change the curve. 
The sequence of $n$ elementary transformations of the first type 
creates one $d$-tacnode of order $n-1$.
Since we apply the second sequence of $n$ elementary transformations 
of the second type on the same fiber, we continue to deepen this
singularity into a $d$-tacnode of order $2n-1$ which is now located
also at the exceptional section.
Hence, when we blow down this section in order to 
return to $\P^2$, this $d$-tacnode of order $2n-1$ 
is blown down too to a more complicated singular point. 
\end{proof} 
 
The description of the curves obtained by Proposition
\ref{gen_2step_coro} is a little bit more complicated: In this case we
indeed have a singular point which is splitted after the first blow-up into
several singular points. 

\begin{prop}\label{sing_gen} 
Let $C$ be a curve of degree $d$. 
Let $(n_1, \cdots, n_k,m_1,\cdots,m_l)\in \N^{k+l}$ be a $(k+l)$-tuple
such that $\sum n_i = \sum m_j$.  
Let $\tilde C = T_{(n_1,\cdots, n_k;m_1,\cdots,m_l)} (C)$ (see Proposition
\ref{gen_2step_coro}). 
 
Then $\tilde C$ has $k+1$ additional singularities to those of $C$:  
$k$ $d$-tacnodes of order $n_i-1$ ($1 \leq i \leq k$), 
and another singular point which  
is a blow down of $l$ $d$-tacnodes of order $m_j-1$ ($1 \leq j \leq l$)  
on the exceptional section of that blow down 
(i.e., the curve has the following 
additional singularities: $[d_{n_1}], \cdots, [d_{n_k}]$  and 
$[d(n_1+\cdots +n_k),([d_{m_1}], \cdots, [d_{m_l}])]$). 
\end{prop} 

\begin{proof} 
The proof is similar to the proofs of the previous propositions.
As before, the elementary transformations of the first type
create $k$ $d$-tacnodes of order $n_i-1$.

The sequences of $m_j$ $(1 \leq j \leq l)$ elementary transformations 
of the second type on the fibers $P_1, \cdots, P_l$ create $l$ 
$d$-tacnodes of order $m_j-1$, which are all located at the 
exceptional section. Hence, when we blow this section down in order to 
return to $\P^2$, these $l$ $d$-tacnodes are 
blown down together to a complicated singular point. 
\end{proof} 

\subsection{Change of the degree of the curve}\label{sec_deg}
In this subsection, we compute the degree of the resulting curve: 

\begin{prop}\label{degree} 
Let $C$ be a plane projective curve of degree $d$.\nl 
Let $(n_1,\cdots,n_k) \in \N^k$. 
Let $\tilde C = T_{(n_1,\cdots, n_k)} (C)$ (see Theorem \ref{gen_2step}).
 
Then the degree of the resulting curve is $d(n_1+\cdots+n_k +1)$.   
\end{prop} 

\begin{proof}
Let $\tilde d$ be the degree of $\tilde C$. 
We have to show that $\tilde d =d(n_1+\cdots+n_k +1)$. 
When we blow up once one of the singularities, say $P$, 
in order to resolve it, we have to decrease the 
self-intersection of the original curve $\tilde C$ by 
$({\rm mult}_{\tilde C} P)^2$ (where ${\rm mult}_{\tilde C} P$ is the local 
multiplicity of $\tilde C$ at $P$) to get the 
self-intersection of the curve $\tilde C$ after the blow-up. 
Since this is the data which is given by the types 
of the singularities, one can compute easily the change 
in the self-intersection. 

So, we start with $\tilde C$ whose
self-intersection is $\tilde d^2$, since $\tilde C$ is in $\P^2$.  
For all $1 \leq i \leq k$, the $n_i$ blow-ups of the singular point 
$[d_{n_i}]$ yield a decreasing of the self-intersection 
by $n_i \cdot d^2$, since the multiplicity of the curve at the
singular point is $d$. The $n_1+\cdots +n_k+1$ blow-ups of the
singular point of the type $[d(n_1+\cdots +n_k),d_{n_1+\cdots +n_k}]$ yield 
an additional decreasing of the self-intersection by  
$(d(n_1+\cdots +n_k))^2+(n_1+\cdots +n_k)d^2$, since the multiplicity 
of the curve at the singular point in the first blow up is 
$d(n_1+\cdots +n_k)$ and in the other blow ups it is again  
$d$. After all these blow-ups, we reach the original curve $C$ 
in $\P^2$ and hence its self-intersection is $d^2$.  
Therefore, we have the following equation: 
$$\tilde d^2 -\sum_{i=1}^k (n_i \cdot d^2) - ((d(n_1+\cdots +n_k))^2+(n_1+\cdots +n_k)d^2)=d^2$$ 
and hence $\tilde d^2 = d^2(n_1+\cdots +n_k+1)^2$, which gives 
us $\tilde d = d(n_1+\cdots +n_k+1)$ as needed. 
\end{proof}
 
One can perform the same computations also for the curves obtained by
the constructions presented in Propositions \ref{gen_2step_coro} 
and \ref{special_case} 
(see Propositions \ref{sing_gen} and \ref{sing_special_case}
respectively for the descriptions of the additional singularities). Hence,
Proposition \ref{degree} holds for those curves too.

\subsection{Families of curves obtained by starting with smooth irreducible curves} 
In this subsection we describe the families of curves which 
are obtained by Uluda\u g's  
original construction and its generalizations if we apply them to a smooth  
irreducible curve $C$ of degree $d$ (and therefore 
$\pi_1(\P^2-C) \cong \Z / d\Z$, see Zariski \cite{Z1}).

\begin{prop}\label{curves1} 
Let $C$ be a smooth irreducible curve of degree $d$. 
Let $n \in \N$. Let $\tilde C$ be the curve obtained by Uluda\u g's
construction in such a way that its  
fundamental group is a central extension of 
$G=\pi _1 (\P ^2 -C)= \Z / d \Z$ by $\Z/(n+1)\Z$.  
 
Then for $n=1$, $\tilde C$ has an intersection point of $d$ smooth 
branches and one $d$-tacnode (i.e., the curve has the singularities:  
$[d]$ and $[d_2]$). 
 
For $n \geq 2$, $\tilde C$ has a $d$-tacnode of order $n-1$,  
and another singular point which is a blow down  
of a $d$-tacnode of order $n-1$ (i.e., the curve has the following 
singularities: $[d_n]$ and $[dn,d_n]$).  
 
The degree of the resulting curve is $d(n+1)$.   
\end{prop} 
 
\begin{proof} 
Since a smooth curve has no singularities, then the only singularities of 
the resulting curve are those   
which were created by  Uluda\u g's construction 
(see Remark \ref{sing_uludag}). 
Therefore the curve has only the singularities described 
in  Remark \ref{sing_uludag}.   

The degree of the curve can be computed directly by Proposition \ref{degree}.    
\end{proof} 
 
For the particular case $d=2$ and $n=1$, we indeed get a quadric 
with a node and a tacnode,  
and its equation can be found in \cite[p. 147, case 2]{Na}: 
$(x^2 +y^2 -3x)^2 = 4x^2(2-x)$. 

Using Proposition \ref{cyclicity}, we have that 
$$\pi _1 (\P ^2 -\tilde C) \cong \Z / (d(n+1))\Z.$$ 
Now, we describe the family of curves which are obtained 
by the general construction (Subsection \ref{second_step}).  
  
\begin{prop} 
Let $C$ be a smooth irreducible curve. 
Let $(n_1,\cdots, n_k) \in \N^k$ be a $k$-tuple.  
Let $\tilde C=T_{(n_1,\cdots, n_k)}(C)$ (see Theorem \ref{gen_2step}).
 
For every $1 \leq i \leq k$, $\tilde C$ has a $d$-tacnode of order $n_i-1$,  
and another singular point which is a blow down  
of a $d$-tacnode of order $n_1+\cdots +n_k-1$ 
(i.e., the curve has the following singularities: 
$[d_{n_1}], \cdots , [d_{n_k}]$ and 
$[d(n_1+\cdots +n_k),d_{n_1+\cdots +n_k}]$).  
 
The degree of the resulting curve is $d(n_1+\cdots +n_k+1)$.   
\end{prop} 
 
\begin{proof} 
Similar to the proof of the previous proposition, but here we use the
results of Proposition \ref{sing_general}.
\end{proof} 
 
Using Proposition \ref{cyclicity} again, we have that 
$$\pi _1 (\P ^2 -\tilde C) \cong \Z / (d(n_1+\cdots +n_k+1))\Z.$$

\subsection{Families of curves obtained by starting with line arrangements} 
In this subsection we describe the families of curves and their groups
which are obtained by Uluda\u g's original construction and 
its generalizations if we apply them to some types of line arrangements.
 
\begin{prop}\label{lines1} 
Let $\cL$ be a line arrangement consisting of $m$ lines intersecting in
one point. Let $n \in \N$. Let $\tilde \cL$ be the curve  
obtained by Uluda\u g's construction in such a way that its 
fundamental group is a central extension of 
$G=\pi _1 (\P ^2 -\cL)= \F_{m-1}$ by $\Z/(n+1)\Z$.  
 
Then for $n=1$, $\tilde \cL$ has two intersection points of $m$ smooth 
branches and one $m$-tacnode (i.e., the curve has the following 
singularities: $[m]$,$[m]$ and $[m_2]$). 
 
For $n \geq 2$, $\tilde \cL$ has one intersection points of $m$ smooth 
branches, one $m$-tacnode of order $n-1$,  
and another singular point which is a blow down  
of a $m$-tacnode of order $n-1$ (i.e., the curve has the following 
singularities: $[m]$,$[m_n]$ and $[mn,m_n]$).  
 
The degree of the resulting curve is $m(n+1)$.   
\end{prop} 
 
\begin{proof} 
Since $\cL$ has one intersection point of $m$ smooth 
branches, then the singularities of 
the resulting curve are those which were created by  
Uluda\u g's construction (see Remark \ref{sing_uludag}) and an
additional singularity which was in $\cL$. 

The degree of the curve is computed directly by Proposition \ref{degree}.    
\end{proof} 
 
Since $H^2(\F_{m-1},\Z/ (n+1)\Z)$ is trivial, then we get that:
$$\pi _1 (\P ^2 -\tilde \cL) \cong  \F_{m-1} \oplus \Z/ (n+1)\Z.$$
Now, we describe the family of curves which are obtained 
by the general construction (Subsection \ref{second_step}).  
  
\begin{prop} 
Let $\cL$ be a line arrangement consisting of $m$ lines intersecting in
one point. Let $(n_1,\cdots, n_k) \in \N^k$ be a $k$-tuple.  
Let $\tilde \cL= T_{(n_1,\cdots, n_k)}(\cL)$ (see Theorem \ref{gen_2step}). 
 
Then: in addition to the original intersection point of $\cL$, 
for every $1 \leq i \leq k$, $\tilde \cL$ has a $m$-tacnode of order $n_i-1$,  
and another singular point which is a blow down  
of a $m$-tacnode of order $n_1+\cdots +n_k-1$ 
(i.e., the curve has the following singularities: 
$[m],[m_{n_1}], \cdots , [m_{n_k}]$ 
and $[m(n_1+\cdots +n_k),m_{n_1+\cdots +n_k}]$).  
 
The degree of the resulting curve is $m(n_1+\cdots +n_k+1)$.   
\end{prop} 
 
\begin{proof} 
Similar to the proof of the previous proposition, but here we use the
results of Proposition \ref{sing_general}.  
\end{proof} 

As before, since $H^2(\F_{m-1},\Z/ (n_1+\cdots +n_k+1)\Z)$ is trivial,
then we get that:
$$\pi _1 (\P ^2 -\tilde \cL) \cong  \F_{m-1} \oplus \Z/ (n_1+\cdots +n_k+1)\Z.$$
Now we deal with another important type of line arrangements: lines in
a general position, which means that there is no intersection of more
than two lines in a point. We describe the family of curves which are obtained 
by the general construction (Subsection \ref{second_step}). 

\begin{prop}\label{general_pos} 
Let $\cL$ be a line arrangement consisting of $m$ lines in a general position.
Let $(n_1,\cdots, n_k) \in \N^k$ be a $k$-tuple.  
Let $\tilde \cL = T_{(n_1,\cdots, n_k)}(\cL)$ (see Theorem \ref{gen_2step}).
 
Then: in addition to the $m \choose 2$ nodal points of $\cL$, 
for every $1 \leq i \leq k$, $\tilde \cL$ has a $m$-tacnode of order $n_i-1$,  
and another singular point which is a blow down  
of a $m$-tacnode of order $n_1+\cdots +n_k-1$ (i.e., the curve has the following 
singularities: $[m_{n_1}], \cdots , [m_{n_k}]$, 
$[m(n_1+\cdots +n_k),m_{n_1+\cdots +n_k}]$ and $m \choose 2$ 
singularities of the type $[2]$).  
 
The degree of the resulting curve is $m(n_1+\cdots +n_k+1)$.   
\end{prop} 
 
\begin{proof} 
Since $\cL$ has $m \choose 2$ nodal points, 
then the singularities of 
the resulting curve are those which were created by  
the general construction (see Proposition \ref{sing_general}) and 
$m \choose 2$ nodal points.

The degree of the curve is computed directly by Proposition \ref{degree}.    
\end{proof} 

Since a central extension of $\Z^{m-1}$ by $\Z/ (n_1+\cdots +n_k+1)\Z$
is not unique, it is interesting to know which group is indeed
obtained in this case. For this, we perform a direct computation for
finding a presentation for the fundamental group of the complement
of $\tilde \cL$, using braid monodromy techniques and van Kampen's theorem
(for similar computations, see \cite{AGT,AmTe,MoTe2}).
We get that: 
$$\pi _1 (\P ^2 -\tilde \cL) \cong \Z^{m-1} \oplus \Z/ (n_1+\cdots +n_k+1)\Z$$
This result is mainly achieved due to the commutative relations induced by
the $m \choose 2$ nodal points, and the torsion subgroup is created by the 
projective relation.

\section{An application to Zariski pairs}\label{zariski_pairs}

As already mentioned, we call {\it a Zariski pair} to a pair of
plane curves which have the 
same combinatorics, but their complements are not homeomorphic. 

In this short section, we want to use the above constructions 
to produce new Zariski pairs.

Not every Zariski pair $(C_1,C_2)$ can produce a family of Zariski pairs
by our construction, since 
even if $G_1=\pi_1(\P^2-C_1)$ and $G_2=\pi_1(\P^2-C_2)$ are
different, it is not guaranteed that $\tilde G_1=\pi_1(\P^2-\tilde C_1)$ 
and $\tilde G_2=\pi_1(\P^2-\tilde C_2)$ will be still different, as
there are several ways to construct the same group by central
extensions. Therefore, we have to characterize Zariski pairs which
induce such families. 

A possible characterization is the following:
\begin{prop}\label{restrict}
Let $(C_1,C_2)$ be a Zariski pair of two irreducible curves.
If $\pi_1(\P^2-C_1)$ is a cyclic group and $\pi_1(\P^2-C_2)$ is not a
cyclic group, then $(\tilde C_1,\tilde C_2)$ is a Zariski pair too. 
\end{prop}

\begin{proof}
Since $(C_1,C_2)$ is a Zariski pair, 
then by definition $C_1$ and $C_2$ have the
same degree and the same singularities. Therefore, using the results
of Section \ref{extension_cyclic}, $\tilde C_1$ and $\tilde C_2$ have the
same degree and the same singularities too.

Since $\pi_1(\P^2-C_1)$ is cyclic and $C_1$ is irreducible, 
then $\pi_1(\P^2-\tilde C_1)$ is cyclic too 
(by Proposition \ref{cyclicity}). On the other hand, a
central extension of non-cyclic group can never be cyclic and hence 
$\pi_1(\P^2-\tilde C_2)$ is not cyclic. 
Therefore, $(\tilde C_1,\tilde C_2)$ is a Zariski pair too.
\end{proof}

The examples of Zariski \cite{Z1,Z2}, Oka \cite{oka2}
and Shimada \cite{shimada2} satisfy the conditions of 
Proposition \ref{restrict}, and hence can be used for producing families of 
new examples of Zariski pairs.

\section*{Acknowledgments} 
I would like to thank Mikhail Zaidenberg and Uzi Vishne for many fruitful 
discussions.  I wish to thank an anonymous referee 
for improving the paper by his suggestions. 
Also, I thank Mikhail Zaidenberg and Institut Fourier for hosting my stay. 
 
The author is partially supported by the Chateaubriand postdoctoral 
fellowship funded by the French government.

\Addresses\recd
 
\end{document}